\documentclass[reqno, 12pt]{amsart}
\usepackage{amsmath}
\usepackage{amssymb}
\usepackage{amsthm}
\usepackage{mathtools}
\usepackage{tikz-cd}
\usepackage{enumerate}
\usepackage[mathscr]{eucal}
\usepackage{graphicx}

\newtheorem{theorem}{Theorem}[section]
\newtheorem{lemma}[theorem]{Lemma}
\newtheorem{proposition}[theorem]{Proposition}
\newtheorem{corollary}[theorem]{Corollary}
\newtheorem{definition}[theorem]{Definition}
\theoremstyle{definition}
\newtheorem{example}[theorem]{Example}
\newtheorem{remark}[theorem]{Remark}

\begin{document}
	\title[]{Uniform convergence of distance functionals under remetrization and infima of hyperspace topologies}
	\author{Yogesh Agarwal and Varun Jindal }
	
	\address{Yogesh Agarwal: Department of Mathematics, Malaviya National Institute of Technology Jaipur, Jaipur-302017, Rajasthan, India}
	\email{yagarwalm247@gmail.com}
	
		\address{Varun Jindal: Department of Mathematics, Malaviya National Institute of Technology Jaipur, Jaipur-302017, Rajasthan, India}
		\email{vjindal.maths@mnit.ac.in}
	\maketitle	
	\begin{abstract}
	The objective of this paper is twofold. In the first half of the paper, we investigate upper parts of the hyperspace convergences determined by uniform convergence of distance functionals on a bornology under different metrizations of a metrizable space. To do this, a new covering property associated with the underlying bornology is introduced. An independent study of this new covering notion in relation to some well-known notions, such as strong uniform continuity, is also presented. In the second half, we study the infima of hyperspace convergences (induced by distance functionals) determined by a family of (uniformly) equivalent metrics. In particular, we establish the existence of the minimum element for the collection of upper Attouch-Wets convergences corresponding to all equivalent metrics on a metrizable space $X$. We show that such a minimum element exists if and only if $X$ has a compatible Heine-Borel metric. Our findings give several new insights into the theory of hyperspace convergences.   
	\end{abstract}
	
	\section{Introduction}
	
	For a metric space $(X,d)$, let $CL(X)$ denote the collection of all nonempty closed subsets of  $(X,d)$. One can identify each $A \in CL(X)$ with the distance functional $d(\cdot, A)$, defined by $d(\cdot,A)(x) = d(x,A) = \inf\{d(x, a):a\in A\}$ for each $x \in X$. Several topologies and convergences have been studied in the literature on $CL(X)$. The metric $d$ on the base space $X$ plays a significant role in constructing several hyperspace topologies. For most of these topologies, the metric space $(X,d)$ can be embedded in $CL(X)$ through the map $x \to \{x\}$. However, in general, the hyperspace $CL(X)$ behaves in a quite different way from space $(X,d)$, that is, two compatible metrics $d, \rho$ on $X$ may not determine the same hyperspace topologies on $CL(X)$. For instance, the \textit{Wijsman convergences} determined by the compatible metrics $d, \rho$ on $X$ need not be identical on $CL(X)$. More generally, if $\mathcal{S}$ is a family of nonempty subsets of $(X,d)$, then the hyperspace convergence on $CL(X)$ defined by $A_{\lambda} \to A \Leftrightarrow d(\cdot, A_{\lambda}) \to d(\cdot, A)$ uniformly on the members of $\mathcal{S}$ is sensitive to the remetrization of $X$. This metric dependence of a hyperspace convergence leads to the following two kinds of problems in the theory of hyperspaces.
	\begin{enumerate}[(A)]
		\item For which compatible metrics $d, \rho$ on $X$ are the corresponding hyperspace convergences equivalent?
		\item To investigate the convergence corresponding to the infima of hyperspace convergences determined by a given family of equivalent metrics on $X$.
	\end{enumerate}  
	This paper aims to investigate the above problems for the upper part of hyperspace convergence determined by uniform convergence of distance functionals on members of a family $\mathcal{S}$. We call this convergence the \textit{$\tau_{\mathcal{S},d}$-convergence} \cite{beer2023bornologies, Idealtopo, cao, cornet1973topologies, sonntagsetconvergences}. Recently, various topological characteristics of $\tau_{\mathcal{S},d}$-convergence and its connections with other hyperspace convergences have been extensively studied by the authors in \cite{agarwal2025set3, agarwal2024set, agarwal2024set2}.

	
	 In the special cases, when $\mathcal{S}$ is the family of all nonempty finite (resp. $d$-bounded) subsets of $X$, the $\tau_{\mathcal{S},d}$-convergence reduces to the classical Wijsman (resp. Attouch-Wets) convergence. These special cases are vital in studying convergences of lower semicontinuous convex functions when identified with their \textit{epigraphs}. Particularly, these are useful in fields such as Variational analysis, Optimization,  Control theory, and other related areas, where the traditional convergences, such as pointwise and uniform convergence of  convex functions, are no longer applicable. For instance, the \textit{Wijsman convergence} for the sequences of convex functions defined on Euclidean spaces is useful for the purposes of stability analysis and well-posedness in optimization \cite{lucchetti2006convexity,wijsconvergence}. In infinite dimensional cases, the \textit{Attouch-Wets convergence} is more helpful in such kinds of problems \cite{attouch1991quantitative,attouch1993}. The classical \textit{Hausdorff distance convergence} naturally appearing in shape analysis \cite{hausshape} and fractal theory \cite{hausfrac} can also be deduced as a special case of $\tau_{\mathcal{S},d}$-convergence for $\mathcal{S}$ to be the family of all nonempty subsets of $X$.
	
	A careful review of the equivalence of classical hyperspace convergences with respect to two compatible metrics on $X$ suggests that for an arbitrary family $\mathcal{S}$, the behavior of  $\tau_{\mathcal{S},d}$-convergence on $CL(X)$  depends on the choice of compatible metric $d$ on $X$ as well as on the nature of family $\mathcal{S}$. Usually such equivalences of hyperspace convergences are described in terms of two inter-related notions: a covering property related to family $\mathcal{S}$; or (strong) uniform continuity of the identity function. Some prominent results and facts in this context are:
	    	    \begin{itemize}
	    	\item The Hausdorff distance convergences corresponding to two equivalent metrics are the same iff the metrics in consideration are uniformly equivalent \cite{ToCCoS}.
	    	
	    	\item The Attouch-Wets convergences corresponding to two equivalent metrics are the same iff the metrics in consideration determine the same class of bounded sets and the same class of bounded real-valued functions on $X$ that are uniformly continuous on  bounded subsets of $X$ \cite{beer1991awucb}.
	    	
	    	\item Equivalently, the Attouch-Wets convergences corresponding to two equivalent metrics $d, \rho$ are the same iff the metrics $d, \rho$ generate the same class of bounded sets and the identity function $id: (X,d) \to (X, \rho)$ is \textit{strongly uniformly continuous} (see, Definition \ref{Stronguniformcontinuity}) on bounded subsets of $X$ and vice-versa \cite{beerattouch}.
	    	
	    	\item Unlike the Attouch-Wets and Hausdorff distance convergences, the equivalence of the Wijsman convergences with respect to two equivalent metrics does not relate naturally to the uniform continuity of the identity function. Even non uniformly equivalent metrics may determine the same Wijsman convergences! (see, Theorem 4.3 of \cite{LechickiWijsman}).
	    	
	    	\item Constatini, Levi, and Zieminska in \cite{costantini1993metrics} characterized the equivalence of Wijsman convergences corresponding to two equivalent metrics on a metrizable space. The resulting condition turns out to be a nice covering condition (\textit{strictly $d$-included}) in terms of the balls of the metrics under consideration.
	    	\item Generalizing the idea of strictly $d$-included Beer et. al in \cite{beer2013gap} introduced a new covering property (\textit{strictly ($\mathcal{S}-d$) included}) corresponding to a family $\mathcal{S}$ of nonempty subsets of a metric space $(X,d)$. Using this covering property, they characterize the equivalence of upper parts of the weak convergences induced by a family of gap functionals corresponding to two different metrics. Similarly, the bornological convergences ($\mathcal{S}$-convergences) with respect to two compatible metrics $d,\rho$ are the same iff the identity function is strongly uniformly continuous on $\mathcal{S}$ in both directions (Theorem 3.13 of \cite{PbciAWc}).
	    	\item In \cite{normsamewijsman}, the author provided another criterion for the equivalence of Wijsman convergences on closed and convex sets in the context of normed linear spaces. The author then used it to give a new proof for Borwein-Fitzpatrick 's theorem \cite{borweinmosco} concerning the existence of the maximum of Wijsman topologies over all equivalent norms. 
	    \end{itemize}      
	    Several researchers have investigated the infima/suprema of hyperspace convergences determined by a family of equivalent metrics on the underlying metric space \cite{beer1992suprema, beer1993weak, costantiniinfimahausdorff, leviinfimum}. For example, the suprema of Wijsman topologies over all equivalent norms is the popular \textit{Joly-slice topology} (Theorem $2.4.5$, \cite{ToCCoS}), and over all equivalent metrics is the classical \textit{Vietoris topology} (Theorem $2.2.5$, \cite{ToCCoS}). Similarly, the set convergence corresponding to the infima of Wijsman topologies over all equivalent metrics is the well-known \textit{Kuratowski-Painelev$\acute{e}$ convergence} \cite{infinmahyperspace}. It is to be noted that for a metrizable space $X$, the infima and suprema of hyperspace convergences when taken over all equivalent metrics become independent of the choice of a compatible metric on $X$.  Such suprema/infima are particularly significant in studying the measurability of multifunctions \cite{beer1987lcf}. Moreover, when a given collection of hyperspace convergences has the minimum element,  certain intrinsic properties of the underlying metrizable space can be characterized \cite{infinmahyperspace,costantini1993metrics, LechickiWijsman}. For example, a metrizable space is locally compact iff the collection of Wijsman convergences over all equivalent metrics has a minimum  \cite{costantini1993metrics}.
	    
	    The organization of this article is as follows. In the 2nd section, we present some basic definitions and notations, and give some preliminary results. Section 3 is devoted to study a new covering property related to a bornology $\mathcal{S}$ on a metric space $(X,d)$. We call this covering property \textit{strict variationally ($\mathcal{S}$-$d$) included} (see, Definition \ref{Stirctvariationallyincluded}). This property is used in Section 4 to study the coincidence of upper parts of the $\tau_{\mathcal{S},d}$ and $\tau_{\mathcal{S},\rho}$-convergences on $CL(X)$ for two compatible metrics $d, \rho$ on a metrizable space $X$. In section 4, we also explore the connection of this new covering notion with the notion of strong uniform continuity. As a byproduct of this, we get a new result for the coincidence of upper bornological convergences with respect to two compatible metrics (Theorem \ref{bornologicalconvergncesuff}). Section 5 delves into the study of infima of hyperspace convergences corresponding to a given family of equivalent metrics on a metrizable space $X$, with special attention on the classical Attouch-Wets and Wijsman convergences.

	  \section{preliminaries}
	  All metric spaces are assumed to have at least two points. For a metric space $(X,d)$, an \textit{ideal} is a family $\mathcal{S}$ of nonempty subsets of $X$ that is closed under finite union and taking subsets of its members. An ideal $\mathcal{S}$ which forms a cover of $X$ is called a \textit{bornology} on $(X,d)$. The concept of bornology is a generalization of the idea of boundedness to topological spaces. Bornologies with some additional assumptions frequently occur in Functional analysis \cite{hogbe}. Examples of some  important metric bornologies include: 
	  \begin{itemize}
	  	\item $\mathcal{F}(X) =$ the family of all nonempty finite subsets of $X$;
	  	\item  $\mathcal{B}_d(X) =$ the family of all nonempty $d$-bounded subsets of $(X,d)$;
	  	\item $\mathcal{K}(X) = $ the family of all nonempty relatively compact subsets of $(X,d)$;
	  	\item $\mathcal{TB}_{d}(X) = $ the family of all nonempty $d$-totally bounded subsets of $(X,d)$;
	  	\item $\mathcal{P}_0(X) = $ the family of all nonempty subsets of $X$.
	  \end{itemize}
	  
	  For $A \subseteq X$, the \textit{$\epsilon$-enlargement of $A$} in $(X,d)$, denoted by $B_d(A, \epsilon)$, is defined as, $B_d(A, \epsilon) = \{x \in X: d(x, A) < \epsilon\}$.
	  Let $$\mathcal{Z}^+ = \{f \in \mathbb{R}^X: f(x)>0~ \forall~ x \in X\}.$$ For a nonempty subset $A$ of $(X,d)$ and $f \in \mathcal{Z}^+$, the \textit{functional enlargement} of $A$ by $f$, denoted by $B_d(A,f)$ is defined as, $$B_d(A,f) = \bigcup_{a \in A}B_d(a,f(a)) = \{x \in X: d(x,a)< f(a) \text{ for some } a \in A\}.$$ The notion of functional enlargement was introduced to find a hit-and-miss type characterization of the $\tau_{\mathcal{S},d}$-convergence (\cite{agarwal2025set3, agarwal2024set}). A functional enlargement of a set $A$ differs in several aspects from its constant enlargement. For example, for any $\epsilon > 0$, $\overline{A} \subseteq B_d(A, \epsilon)$ but it may not be always true that $\overline{A} \subseteq B_d(A, f)$ for $f \in \mathcal{Z}^+$. Note that for any $f,g \in \mathcal{Z}^+$ and $A \in \mathcal{P}_0(X)$ if $\inf_{x \in A}(g(x)-f(x))>0$, then $\overline{B_d(A,f)}\subseteq B_d(A,g)$. We denote, the collection of all proper functional enlargements of members of a bornology $\mathcal{S}$ on $(X,d)$ by $\mathcal{S}^d_{\mathcal{Z}^+}$. 
	  
	  Consider a bornology $\mathcal{S}$ on $(X,d)$ and $A \in \mathcal{P}_0(X)$. Then $A$ is called \textit{weakly $\mathcal{S}$-totally bounded} (resp. \textit{$\mathcal{S}$-totally bounded}) if for any $\epsilon >0$, $\exists$ $S \in \mathcal{S}$ (resp. $S \in \mathcal{S}$ and $S \subseteq A$) satisfying $A \subseteq B_d(S,\epsilon)$.  
	The collection of all weakly $\mathcal{S}$-totally bounded (resp. $\mathcal{S}$-totally bounded) sets is denoted by $\mathcal{S}^*$(resp. $\mathcal{S}_*$) \cite{beer2009totalboundedness, BCL}. To specify the metric, we use the notations $\mathcal{S}_d^*$ and $\mathcal{S}^d_*$ whenever required.
	

A metric space $(X,d)$ is said to be an \textit{almost convex metric space (a.c.m.s.)} if for each $x \in X$ and for each $\alpha, \epsilon > 0$, we have $B_d(B_d(x,\alpha), \epsilon) = B_d(x, \alpha+ \epsilon)$ \cite{ToCCoS}. A prime example of an a.c.m.s. is a normed linear space.

Note that for $A \in \mathcal{P}_0(X)$ and $f \in \mathcal{Z}^+$ and $\epsilon > 0$, we always have 
$$B_d(B_d(A,f), \epsilon) \subseteq B_d(A, f+\epsilon).$$It can be easily seen that in an a.c.m.s., we have $B_d(B_d(A,f), \epsilon) = B_d(A, f+\epsilon)$.

	  \begin{definition}\normalfont(\cite{Idealtopo,cao,cornet1973topologies})
	  	Let $(X,d)$ be a metric space and $\mathcal{S} \subseteq \mathcal{P}_0(X)$. The topology $\tau_{\mathcal{S},d}^+$ on $\mathcal{P}_0(X)$ has a neighborhood base at any $A \in \mathcal{P}_0(X)$ consisting of sets of the form $$ [S,\epsilon]^+(A) = \{C \in \mathcal{P}_0(X) :  d(x,A) - d(x,C) < \epsilon ~\forall x \in S\},$$ where $S \in \mathcal{S}$ and $\epsilon >0$.
	  So a net $(A_\lambda)$ is $\tau_{\mathcal{S},d}^+$-convergent to $A$ if and only if for any $S \in \mathcal{S}$ and $\epsilon > 0$, eventually, $ d(x,A) - d(x, A_\lambda) < \epsilon$ for all $x \in S$.
	  	\end{definition}	
	  	 	 Similarly, the topology $\tau_{\mathcal{S},d}^-$ on $\mathcal{P}_0(X)$ has a neighborhood base at any $A \in \mathcal{P}_0(X)$ having sets of the form 
	  	 	$$ [S,\epsilon]^-(A) = \{C \in \mathcal{P}_0(X) : d(x,C) - d(x,A) < \epsilon ~~ \forall x \in S \}.$$

	  The topology $\tau_{\mathcal{S},d}$ is the supremum of $\tau_{\mathcal{S},d}^-$ and $\tau_{\mathcal{S},d}^+$. It is a uniformizable topology for which the family $\{U_{\mathcal{S},\epsilon} : S \in \mathcal{S}, \epsilon > 0\}$ forms a base for a compatible uniformity, where $$ U_{S,\epsilon} = \{(A,C) \in \mathcal{P}_0(X) \times \mathcal{P}_0(X) : |d(x,A) - d(x,C)| < \epsilon\ \ \forall x \in S\}.$$

	  
	  An important representation for $\tau_{\mathcal{S},d}^+$-convergence was proved in \cite{agarwal2024set}. Through this representation, one can view $\tau_{\mathcal{S},d}^+$-convergence as a miss-type convergence using functional enlargements of members of $\mathcal{S}$. This viewpoint is quite useful in our investigations. For the readers' convenience, we record its statement below.   
	  \begin{theorem}\label{Miss theorem for upper Tsd} \normalfont{(Theorem 4.1 of \cite{agarwal2024set})}
	  	Let $(X,d)$ be a metric space and let $\mathcal{S}$ be a bornology on $X$. Suppose $(A_\lambda)$ is a net in $CL(X)$ and $A \in CL(X)$. Then the following statements are equivalent: 
	  	\begin{enumerate}[(i)]
	  		\item for $S \in \mathcal{S}$ and $f, g \in \mathcal{Z}^+$ with $\inf \{g(x) - f(x) : x \in S\} > 0$, whenever $A$ misses $g$-enlargement of $S$, then $(A_\lambda)$ misses the $f$-enlargement of $S$ eventually; 
	  		
	  		\item $(A_\lambda)$ is $\tau_{\mathcal{S},d}^+$-convergent to $A$. 
	  	\end{enumerate}
	  \end{theorem}
	  
	  The Attouch-Wets convergence and the Hausdorff distance convergence can also be deduced as particular cases of another set convergence, called \textit{bornological convergence}, introduced by Lechicki, Levi, and Spakowski (\cite{BCL}) and further studied in \cite{beer2023bornologies,PbciAWc,Idealtopo}.  
	  \begin{definition}\normalfont(\cite{BCL})
	  	Suppose $\mathcal{S}$ is an ideal of subsets of a metric space $(X,d)$. Suppose $(A_\lambda)$ is a net in $\mathcal{P}_0(X)$ and $A \in \mathcal{P}_0(X)$.
	  	
	  	\begin{enumerate}[(i)]
	  		\item The net $(A_\lambda)$ is said to be \textit{lower bornological convergent} to $A$ in $\mathcal{P}_0(X)$ denoted by, $A_\lambda$ is $(\mathcal{S}$-$d)^-$-convergent to $A$ if for each $\epsilon > 0$ and $S \in \mathcal{S}$, eventually, $ A \cap S \subseteq B_d(A_\lambda, \epsilon)$. 
	  	\item The net $(A_\lambda)$  is said to be \textit{upper bornological convergent} to $A$ in $\mathcal{P}_0(X)$ denoted by, $A_\lambda$ is $(\mathcal{S}$-$d)^+$-convergent to $A$ if for each $\epsilon > 0$ and $S \in \mathcal{S}$, eventually, $ A_\lambda \cap S \subseteq B_d(A, \epsilon)$. 
	  	\end{enumerate}

	  Finally, the net $(A_\lambda)$ is said to be $(\mathcal{S}$-$d)$-convergent to $A$ in $\mathcal{P}_0(X)$ if for each $\epsilon > 0$ and $S \in \mathcal{S}$, eventually, we have $A_\lambda \cap S \subseteq B_d(A, \epsilon)\ \ \text{and}\ \ A \cap S \subseteq B_d(A_\lambda,\epsilon)$.
	  \end{definition}
	  	 
	  	 The following notion of uniform continuity was formally introduced in \cite{Suc}. However, this notion was used much earlier by Beer and Concilio in \cite{beer1991awucb} to characterize the equivalence of Attouch-Wets topologies corresponding to different compatible metrics on a metrizable space.   
	  \begin{definition}\normalfont\label{Stronguniformcontinuity}
	  	Let $(X,d)$ and $(Y, \rho)$ be two metric spaces. Suppose $f:X \to Y$ is a function and $S \in \mathcal{P}_0(X)$. Then $f$ is said to be \textit{strongly uniformly continuous} on $S$ if for each $\epsilon >0$ there exists a $\delta > 0$ such that whenever $\{x,y\}\cap S \neq \emptyset$ with $d(x,y) < \delta$, then $\rho(f(x), f(y)) < \epsilon$.   
	  \end{definition}
	    If the function $f$ is strongly uniformly continuous on each $S \in \mathcal{S}$, then we say $f$ is strongly uniformly continuous on $\mathcal{S}$.
	 \section{Covering properties}
	    Covering properties naturally arise in the theory of hyperspaces while investigating the equivalence of a given hyperspace convergence with respect to two compatible metrics. In this section, we introduce a new covering property corresponding to a bornology $\mathcal{S}$ on $(X,d)$. This new notion differs from the existing ones in the sense that, in its definition, we require functional enlargements of members of $\mathcal{S}$, whereas the existing covering properties use the constant enlargements. In the present situation, the motivation of using such enlargements comes from the miss-type representation of $\tau_{\mathcal{S},d}^+$-convergence (Theorem \ref{Miss theorem for upper Tsd}). The equivalence of this new covering property with the existing covering properties allows us to view $\tau_{\mathcal{S},d}^+$-convergence as a convergence corresponding to the weak topology generated by certain geometric set functionals.    
	    
%
	   
	   We first recall definitions of existing covering properties.
	   \begin{definition}\normalfont(\cite{costantini1993metrics})\label{Strictlyinclued}
	   	Let $(X,d)$ be a metric space. Let $A,C$ be two nonempty  subsets of $X$. Then $A$ is said to be \textit{strictly $d$-included} in $C$ if there exist $x_1, \ldots, x_n \in X$ and $0 < \alpha_i < \epsilon_i$, where $i = 1, \ldots, n$ such that $$A \subseteq \cup_{i=1}^nB_d(x_i, \alpha_i) \subseteq \cup_{i=1}^nB_d(x_i, \epsilon_i) \subseteq C.$$	   \end{definition}
	   	
	   	The notion of strictly $d$-included was introduced by C. Costantini et al. in \cite{costantini1993metrics} to characterize the equivalence of the Wijsman topologies corresponding to two equivalent metrics $d$ and $\rho$ on a metrizable space $X$.

\begin{definition}\normalfont(\cite{beer2013gap}) \label{StrictlySdinclued}
	 Let $(X,d)$ be a metric space and let $\mathcal{S}$ be a bornology on $(X,d)$. Suppose $A$ and $C$ are two nonempty subsets of $X$ such that $A \subseteq C$. Then $A$ is said to be \textit{strictly $(\mathcal{S}$-$d)$ included} in $C$ if  there exist an $S_1, \ldots, S_n \in \mathcal{S}$ and $0< \alpha_i < \epsilon_i$, $i = 1, \ldots,n$ satisfying
	 $$A \subseteq \cup_{i=1}^nB_d(S_i, \alpha_i) \subseteq \cup_{i=1}^nB_d(S_i,\epsilon_i) \subseteq C.$$
\end{definition}
	The notion of strictly $(\mathcal{S}$-$d)$ included was introduced by G. Beer et al. in \cite{beer2013gap} to study the equivalence of \textit{gap topologies} corresponding to different metrics. This notion also turned out to be helpful in studying the weak representation of $\tau_{\mathcal{S},d}^+$-convergence by the family $\{D_d(S, \cdot): S	 \in \mathcal{S}\}$ of gap functionals \cite{agarwal2024set}. 
	
	The following new notion is important in our investigation.  
	\begin{definition}[Strict Variationally $(\mathcal{S}$-$d)$ inclusion]\normalfont \label{Stirctvariationallyincluded}
  Let $(X,d)$ be a metric space and let $\mathcal{S}$ be a bornology on $(X,d)$. Suppose $A$ and $C$ are two nonempty subsets of $X$ such that $A \subseteq C$. Then $A$ is said to be \textit{strict variationally $(\mathcal{S}$-$d)$ included} in $C$ if  there exist an $S \in \mathcal{S}$ and $f,g \in \mathcal{Z}^+$ with $\inf_{x \in S}(g(x)-f(x))>0$ such that
  $$A \subseteq B_d(S, f) \subseteq B_d(S,g) \subseteq C.$$
	\end{definition}
	In other words, if $A$ can be covered by a functional enlargement of a member of $\mathcal{S}$ and further if that enlargement is fattened by some small margin, it still remains in $C$. Notice that existence of such $S \in \mathcal{S}$, and $f,g \in \mathcal{Z}^+$ depends on both $A$ and $C$.   
	
	Clearly, in a metric space $(X,d)$, 
	
	\[
	\begin{tikzcd}
		A \text{ is strictly $d$-included in } C \arrow[r,Rightarrow]   & A \text{ is strictly } (\mathcal{S}\text{-}d) \text{ included in } C \arrow[d, Rightarrow] \\
		 &\hspace{-1cm} A \text{ is strict variationally } (\mathcal{S}\text{-}d) \text{ included in } C\end{tikzcd}
	\]

 However, the reverse implications need not be true in general. We show it by the following example. 
	
	\begin{example}
		Let $(X,d) = (\mathbb{R}^2, d_e)$, where $d_e$ represents the Euclidean metric and let $\mathcal{B} = \{B \subseteq \mathbb{R}^2 : \text{ there exists } M > 0 \text{ such that } \vert y \vert \leq M \text{ for all } (x,y) \in B\}$. Clearly, $\mathcal{B}$ is a bornology on $\mathbb{R}^2$. Let $A$ be the region $\{(x,y) \in \mathbb{R}^2 : \vert y \vert \leq \frac{x}{3}, x \geq 5\}$ and $C$ be the region $\{(x,y) \in \mathbb{R}^2: y \leq 3 x, x \geq 0\}$ in Euclidean plane. We show that $A$ is strict variationally $(\mathcal{B}$-$d)$ included in $C$ but not strictly $(\mathcal{B}$-$d)$ included in $C$. Assume $B_0 = \{(x,0) : x \geq 5\}$, clearly $B_0 \in \mathcal{B}$. Take $f \in \mathcal{Z}^+$ such that $f(x,y) = \begin{cases}
			\frac{ x - 1}{2} & x > 1\\
			1                   & \text{otherwise}
		\end{cases}$. Note that $B_d(B_0,f) \neq X$ as $\{(x,0) : x < 0\} \cap  B_d(B_0,f) = \emptyset$. Consider $(x,y) \in A$, then $d_e((x,0), (x,y)) = \vert y \vert \leq \frac{x}{3} < \frac{x}{2}-\frac{1}{2}$ as $x \geq 5$. So $A \subseteq B_d(B_0, f)$. Now, define $g \in \mathcal{Z}^+$ by $g(x,y) = \begin{cases}
			\frac{ x }{2}   & x \geq 1\\
			1                   & \text{otherwise}
		\end{cases}  $. 
		Clearly, $\inf\{g(x,y) - f(x,y) : (x,y) \in B_0\}  = \frac{1}{2} > 0$. Take $(x_0, y_0) \in B_d(B_0,g)$. Then for some $(x,0) \in B_0$, $d_e((x_0,y_0) , (x,0)) < g(x,0)$. So $\sqrt{(x-x_0)^2+y_0^2} < \frac{x}{2}$. Then $\frac{x}{2} < x_0 < \frac{3x}{2}$ and $ -\frac{x}{2} < y_0 < \frac{x}{2}$. Consequently,  $y_0 < \frac{x}{2} < x_0 < 3 x_0$ thus, $(x_0, y_0) \in C$. This gives, $B_d(B_0,g) \subseteq C$. Hence, $A$ is strict variationally $(\mathcal{B}$-$d)$ included in $C$. Now consider any $B \in \mathcal{B}$ and $\epsilon > 0$, then there is an $M > 1$ such that $\vert y \vert \leq M$. Take $x > 5(M + 2 \epsilon)$ and $y = M + 2 \epsilon$. Then $(x,y) \in A$ but not belong to $B_d(B,\epsilon)$. So $A$ is not strictly $(\mathcal{B}$-$d)$ included in $C$.\qed 
			\end{example}
		
		Note that when $\mathcal{S} = \mathcal{F}(X)$, all three covering properties become equivalent. We now examine the coincidence of the latter two covering properties. It turns out that this coincidence gives us a weak representation of $\tau_{\mathcal{S},d}^+$-convergence by the family $\{D_d(S, \cdot): S \in \mathcal{S}\}$ of gap functionals.  	
		Recall that for a nonempty subset $A$ of $(X,d)$, the \textit{gap functional corresponding to $A$}, denoted by $D_d(A, \cdot)$ is a real-valued function on $\mathcal{P}_0(X)$ defined by, $D_d(A, \cdot)(B)= D_d(A, B)= \inf\{d(a,B): a \in A\}$, where $B \in \mathcal{P}_0(X)$.
			\begin{theorem}
				Let $(X,d)$ be a metric space and let $\mathcal{S}$ be a bornology on $X$. Then the following statements are equivalent:
				\begin{enumerate}[(i)]
					\item $\tau_{\mathcal{S},d}^+$ is the weakest topology on $CL(X)$ such that each member of the family $\{D_d(S, \cdot): S \in \mathcal{S}\}$ is lower semi-continuous;
					\item for each pair $(A,C)$ of nonempty proper subsets of $X$ such that $A\subseteq C$, $A$ is strict variationally $(\mathcal{S}$-$d)$ included in $C$ $\Leftrightarrow$ $A$ is strictly $(\mathcal{S}$-$d)$ included in $C$.
				\end{enumerate}
			\end{theorem}
			\begin{proof}
				$(i)\Rightarrow (ii)$. Suppose $A,C$ are two nonempty proper subsets of $X$ such that $A$ is strict variationally $(\mathcal{S}$-$d)$ included in $C$. Choose $S \in \mathcal{S}$ and $f,g \in \mathcal{Z}^+$ with $\inf_{x \in S}(g(x)-f(x))>0$ satisfying, $$A \subseteq B_d(S,f) \subseteq B_d(S,g)\subseteq C.$$ By Theorem $5.5$ of \cite{agarwal2024set}, there exist $S_i \in \mathcal{S}$ and $0 < \alpha_i < \epsilon_i$, $i = 1, \ldots, n$ satisfying $$B_d(S,f) \subseteq \cup_{i=1}^nB_d(S_i, \alpha_i) \subseteq \cup_{i=1}^nB_d(S_i, \epsilon_i)\subseteq B_d(S,g).$$ Consequently, $A$ is strictly $(\mathcal{S}$-$d)$ included in $C$.
				
				$(ii) \Rightarrow (i)$. By Theorem $5.5$ of \cite{agarwal2024set} and $(ii)$, it is enough to show that for each $S \in \mathcal{S}$ and $f,g \in \mathcal{Z}^+$ such that $\inf_{x \in S}(g(x)-f(x)) =r>0$ and $B_d(S,g) \neq X$, $B_d(S,f)$ is strict variationally $(\mathcal{S}$-$d)$ included in $B_d(S,g)$. Define $p = f+ \frac{r}{4}$ and $q = f + \frac{r}{2}$. Then $p,q \in \mathcal{Z}^+$ and $$B_d(S,f) \subseteq B_d(S,p) \subseteq B_d(S,q) \subseteq B_d(S,g).$$  
 			\end{proof}
 		
 		The equivalence of strictly $d$-inclusion and strict variationally ($\mathcal{S}$-$d$) inclusion yields the weak representation of $\tau_{\mathcal{S},d}^+$-convergence on $CL(X)$ by the family $\{d(x,\cdot): x \in X\}$ of distance functionals. The proof of this depends on a later result (Theorem \ref{Tsddiffbornologiesmetrics}).

			We end this section, by exploring the relation of strict variationally $(\mathcal{S}$-$d)$ inclusion with the notion of $\mathcal{S}$-totally bounded sets.
			\begin{proposition}\label{SVIP1}
				Let $(X,d)$ be a metric space and let $\mathcal{S}$ be a bornology on $X$. Suppose $A,C \in \mathcal{P}_0(X)$. Then the following statements are true.
				\begin{enumerate}[(i)]
					\item If $A$ is $\mathcal{S}$-totally bounded, then for any $\epsilon >0$ $A$ is strict variationally $(\mathcal{S}$-$d)$ included in $B_d(A, \epsilon)$;
					\item If $A$ is strict variationally $(\mathcal{S}$-$d)$ included in $C$, then $D_d(A, X\setminus C) > 0$. The converse holds if $A$ is $\mathcal{S}$-totally bounded. 
			\end{enumerate}
				\end{proposition}
				
			\begin{proof}$(i)$. Let $\epsilon > 0$. Then there exists $S \in \mathcal{S}$ such that $S \subseteq A$ and $A \subseteq B_d(S, \frac{\epsilon}{3})$. So we have $$A \subseteq B_d(S, \frac{\epsilon}{3})\subseteq B_d(S, \frac{\epsilon}{2}) \subseteq B_d(A, \epsilon).$$
				
				$(ii)$. Suppose $A$ is strict variationally $(\mathcal{S}$-$d)$ included in $C$. Then $\exists$ $S \in \mathcal{S}$ and $f,g \in \mathcal{Z}^+$ such that $\inf_{x \in S}(g(x)-f(x))= \epsilon>0$ and $$A \subseteq B_d(S,f) \subseteq B_d(S,g) \subseteq C.$$ So $B_d(A, \epsilon) \subseteq C$. Consequently, $D_d(A, X\setminus C) > 0$. 
				
				Conversely, suppose $A$ is $\mathcal{S}$-totally bounded and $D_d(A, X\setminus C) > 0$. Let $r = D_d(A, X\setminus C)$. By $(i)$, $A$ is strict variationally $(\mathcal{S}$-$d)$ included in $B_d(A,r)$. Thus, $A$ is strict variationally $(\mathcal{S}$-$d)$ included in $C$.  			
			\end{proof}
			\section{Uniform convergence of distance functionals under remetrization}
		This section is devoted to study the equivalence of $\tau_{\mathcal{S},d}^+$ and $\tau_{\mathcal{S},\rho}^+$-convergences corresponding to compatible metrics $d, \rho$ on a metrizable space $X$. The notions of strict variationally $(\mathcal{S}$-$d)$ inclusion and of functional enlargements of members of $\mathcal{S}$ play a vital role in our investigation. Further, motivated by the condition for the equivalence of the Attouch-Wets convergences with respect to two compatible metrics, we relate the equivalence of $\tau_{\mathcal{S},d}^+$ and $\tau_{\mathcal{S},\rho}^+$-convergences with the notion of strong uniform continuity. In fact, a necessary condition for the inclusion $\tau_{\mathcal{S},d}^+ \supseteq \tau_{\mathcal{S},\rho}^+$ comes out to be the strong uniform continuity of the identity function $Id: (X,d) \to (X, \rho)$ on $\mathcal{S}$. Consequently, by Theorem 3.13 of \cite{PbciAWc},  $\tau_{\mathcal{S},d}^+ \supseteq \tau_{\mathcal{S},\rho}^+$ implies $(\mathcal{S}-d)^+ \geqslant (\mathcal{S}-\rho)^+$ on $CL(X)$. Also the reverse implication holds for the bornologies which are stable under proper functional enlargements. Interestingly, under certain convexity assumptions on the compatible metrics $d, \rho$ on $X$, we obtain that the strong uniform continuity of the identity function from $(X,d)$ to $(X, \rho)$ on the proper functional enlargements of members of $\mathcal{S}$ is sufficient for the equivalence of $\tau_{\mathcal{S},d}^+$ and $\tau_{\mathcal{S},\rho}^+$-convergences.      	
			\begin{theorem}\label{upperTsdcompatiblemetric}
				Let $X$ be a metrizable space and let $\mathcal{S}$ be a bornology on $X$. Suppose $d, \rho$ are two compatible metrics on $X$. Then the following statements are equivalent:
				\begin{enumerate}[(i)]
					\item $\tau_{\mathcal{S},d}^+ \supseteq \tau_{\mathcal{S},\rho}^+$ on $CL(X)$;
					\item  for each $S \in \mathcal{S}$ and $f, g \in \mathcal{Z}^+$ such that  $\inf\{g(x) - f(x) : x \in S\} > 0$ and $B_\rho(S,g) \neq X$, then $B_\rho(S, f)$ is strict variationally $(\mathcal{S}$-$d)$ included in $B_\rho(S, g)$;
					\item for each pair $(A,C)$ of nonempty proper subsets of $X$ such that $A$ is strict variationally $(\mathcal{S}$-$\rho)$ included in $C$, we have $A$ is strict variationally $(\mathcal{S}$-$d)$ included in $C$.     
					\end{enumerate} 
			\end{theorem}
	\begin{proof}
			$(i) \Rightarrow (ii)$.
			Suppose $(ii)$ fails. Then there exists $S_0 \in \mathcal{S}$ and $f ,g \in \mathcal{Z}^+$ such that $\inf\{g(x) - f(x): x \in S_0\} > 0$ and $B_\rho(S_0, g) \neq X$ but $B_\rho(S_0,f)$ is not strict variationally $(\mathcal{S}$-$d)$ included in $B_\rho(S_0,g)$. So for each $S \in \mathcal{S}$ and $p,q \in \mathcal{Z}^+$ such that $\inf_{x \in S}(q(x)-p(x))>0$ and $B_d(S,q) \subseteq B_{\rho}(S_0,g)$, there exists $y_{S,p}\in B_{\rho}(S_0,f)\setminus B_d(S,p)$.
		 
			Define \begin{equation*}
				\Sigma = \left\{(S, p) \in \mathcal{S} \times \mathcal{Z}^+ :\begin{split} \exists q \in \mathcal{Z}^+ \text{ with } \inf_{x \in S}(q(x) - p(x)) > 0 \\ \text{and } B_d(S, q) \subseteq B_\rho(S_0, g)\end{split}\right\}.\end{equation*}
			To see $\Sigma$ is nonempty, choose $x_0 \in S_0$. Since $d,\rho$ are compatible metrics, there exists $0< \alpha_{x_0} < \epsilon_{x_0} $ such that $B_d(x_0, \epsilon_{x_0}) \subseteq B_\rho(x_0, g(x_0))$. Define $\alpha \in \mathcal{Z}^+$ such that $\alpha(x_0) = \alpha_{x_0}$ and $\alpha(x) = 1$ otherwise. Then $(\{x_0\}, \alpha) \in \Sigma$.   Direct $\Sigma$ as follows: $(S_1, p_1) \leqslant (S_2, p_2)$ if $S_1 \subseteq S_2$ and $p_1(x) \leq p_2(x)$ $\forall x \in S_1$. It is not difficult to verify that $(\Sigma, \leqslant)$ is a directed set. Define $A = X \setminus B_\rho(S_0, g)$ and for every $(S, p) \in \Sigma$, $A_{(S, p)} = A \cup \{y_{S,p}\}$, where $y_{S,p}\in  B_{\rho}(S_0,f)\setminus B_d(S,p)$. Clearly, $A \in CL(X)$ and $(A_{(S, p)})_{(S, p) \in \Sigma}$ is a net in $CL(X)$. We claim that the net $(A_{(S, p)})$ is $\tau_{\mathcal{S},d}^+$-convergent to $A$ but it is not $\tau_{\mathcal{S},\rho}^+$-convergent to $A$. Take any $S' \in \mathcal{S}$ and $m, l\in \mathcal{Z}^+$ with $\inf\{l(x) - m(x) : x \in S' \} > 0$ and $A \cap B_d(S',l)  = \emptyset$. So $(S',m) \in \Sigma$. Then for $(S, p) \in \Sigma$ with $(S, p) \geqslant (S', m)$, we have $A_{(S, p)} \cap B_d(S', m) = \emptyset$. Consequently, by Theorem \ref{Miss theorem for upper Tsd}, the net $(A_{(S, p)})$ is $\tau_{\mathcal{S},d}^+$-convergent to $A$.

			 Now, $A \cap B_\rho(S_0, g)  = \emptyset$ but for any $(S, p) \in \Sigma$, we have $A_{(S, p)} \cap B_\rho(S_0, f)  \neq \emptyset$. Hence, by Theorem \ref{Miss theorem for upper Tsd}, the net $(A_{(S, p)})$ is not $\tau_{\mathcal{S},\rho}^+$-convergent to $A$. This is a contradiction to our hypothesis. 
			
			$(ii) \Rightarrow(i)$.  Suppose $(A_\lambda)$ is a net which is $\tau_{\mathcal{S},d}^+$-convergent to $A$ in $CL(X)$. We need to show that the net $(A_\lambda)$ is $\tau_{\mathcal{S},\rho}^+$-convergent to $A$ in $CL(X)$. If $A = X$, then we are done. Otherwise, take $S \in \mathcal{S}$ and $g \in \mathcal{Z}^+$ such that $A \cap B_\rho(S, g) = \emptyset$. By $(ii)$, for any $f \in \mathcal{Z}^+$ such that $\inf_{x \in S}(g(x) - f(x)) > 0$, there is an $S' \in \mathcal{S}$ and $p,q \in \mathcal{Z}^+$ with $\inf_{x \in S'}(q(x) - p(x)) > 0$ satisfying $$B_\rho(S, f) \subseteq B_d(S',p) \subseteq B_d(S', q) \subseteq B_\rho(S, g).$$ 
			So $A \cap B_d(S', q) = \emptyset$. By the hypothesis, we have $A_\lambda \cap B_d(S', p) = \emptyset$ eventually. Consequently,  $A_\lambda \cap B_\rho(S, f) = \emptyset$ eventually. Hence, the net $(A_\lambda)$ is $\tau_{\mathcal{S},\rho}^+$-convergent to $A$.
			
			The implication $(ii) \Rightarrow (iii)$ is easy to verify and $(iii) \Rightarrow (ii)$ can be proved by taking $A = B_{\rho}(S,f)$ and $C = B_{\rho}(S,g)$.
			\end{proof}
			
			\begin{corollary}
				Let $d, \rho$ be two compatible metrics on a metrizable space $X$ and let $\mathcal{S}$ be a bornology on $X$. Then the following statements are equivalent:
				\begin{enumerate}[(i)]
					\item $\tau_{\mathcal{S},d}^+ = \tau_{\mathcal{S},\rho}^+$ on $CL(X)$;
					
					\item for each pair $(A,C)$ of nonempty proper subsets of $X$, $A$ is strict variationally $(\mathcal{S}$-$\rho)$ included in $C$ $\Leftrightarrow$ $A$ is strict variationally $(\mathcal{S}$-$d)$ included in $C$.     
				\end{enumerate}  
			\end{corollary}
		The next corollary characterizes $\tau_{\mathcal{S},d}^+ \supseteq \tau_{\mathcal{S},\rho}^+$ in terms of the gap between functional enlargements of members of $\mathcal{S}$. 	
	\begin{corollary}
		Let $X$ be a metrizable space and let $\mathcal{S}$ be a bornology on $X$. Suppose $d,\rho$ are two equivalent metrics on $X$. If each proper $\rho$-functional enlargement of members of $\mathcal{S}$ is $\mathcal{S}$-totally bounded with respect to $d$, then the following statements are equivalent:
		\begin{enumerate}
			\item [(i)] $\tau_{\mathcal{S},d}^+ \supseteq \tau_{\mathcal{S},\rho}^+$ on $CL(X)$;
			\item [(ii)] for each $S \in \mathcal{S}$ and $f,g \in \mathcal{Z}^+$ whenever $\inf_{x \in S}(g(x)-f(x))>0$ and $B_{\rho}(S,g) \neq X$, then $D_d (B_\rho(S,f), X\setminus(B_\rho(S,g))>0$.    
		\end{enumerate}
	\end{corollary}
	\begin{proof}
		$(i) \Rightarrow (ii)$. Suppose $\tau_{\mathcal{S},d}^+ \supseteq \tau_{\mathcal{S},\rho}^+$. Then by Theorem \ref{upperTsdcompatiblemetric}, for each $S \in \mathcal{S}$ and $f,g \in \mathcal{Z}^+$ such that $\inf_{x \in S}(g(x)-f(x))>0$ and $B_\rho(S,g)\neq X$, $B_\rho(S,f)$ is strict variationally $(\mathcal{S}$-$d)$ included in $B_\rho(S,g)$. By Proposition \ref{SVIP1} $(ii)$, $D_d (B_\rho(S,f), X\setminus(B_\rho(S,g)) > 0$.
		
		$(ii) \Rightarrow (i)$. Let $S \in \mathcal{S}$ and $f,g \in \mathcal{Z}^+$ such that $\inf_{x \in S}(g(x)-f(x))> 0$ and $B_\rho(S,g)\neq X$. Since $B_\rho(S,f)$ is $\mathcal{S}$-totally bounded with respect to $d$, by Proposition \ref{SVIP1} $(ii)$, $B_\rho(S,f)$ is strict variationally $(\mathcal{S}$-$d)$ included in $B_\rho(S,g)$. Consequently, by Theorem \ref{upperTsdcompatiblemetric}, $\tau_{\mathcal{S},d}^+ \supseteq \tau_{\mathcal{S},\rho}^+$ on $CL(X)$.   
	\end{proof}
	 We now deduce particular cases using Theorem \ref{upperTsdcompatiblemetric} and results of \cite{ToCCoS}. 
	\begin{corollary}
		Let $X$ be a metrizable space and let $d, \rho$ be compatible metrics on $X$. Then the following are equivalent:
		 \begin{enumerate}[(i)]
		 	\item  $\tau_{W_d} \supseteq \tau_{W_\rho}$ on $ CL(X)$;
		 	\item $\tau_{W_d}^+ \supseteq \tau_{W_\rho}^+$ on $ CL(X)$;
		 	\item for each $x \in X$ and $ 0 < \alpha < \epsilon$ such that $B_\rho(x, \epsilon) \neq X$, then $B_\rho(x, \alpha)$ is strict variationally $(\mathcal{F}(X)$-$d)$ included in $B_\rho(x, \epsilon)$;
		 	\item for each $x \in X$ and $ 0 < \alpha < \epsilon$ such that $B_\rho(x, \epsilon) \neq X$, then $B_\rho(x, \alpha)$ is strictly $d$-included in $B_\rho(x, \epsilon)$.
		 \end{enumerate}
	\end{corollary}
	\begin{proof}
		Since the lower Wijsman convergence is metric-independent, the equivalence $(i) \Leftrightarrow (ii)$ holds. The equivalence $(ii) \Leftrightarrow (iii)$ follows from Theorem \ref{upperTsdcompatiblemetric}. Finally, $(iii) \Leftrightarrow (iv)$ follows from Definitions \ref{Strictlyinclued} and \ref{Stirctvariationallyincluded}. 
	\end{proof}
	In order to deduce equivalence of Hausdorff distance convergence with respect to different compatible metrics, we first recall Efremovic Lemma. The Efremovic Lemma states that if $(x_n)$ and $(y_n)$ are two sequences in a metric space $(X,d)$ such that for each $n$, $d(x_n,y_n) > \epsilon$, where $\epsilon > 0$,  then there exist subsequences $(x_{n_k})$ and $(y_{n_k})$ satisfying $D_d(\{x_{n_k}: k\in \mathbb{N}\}, \{y_{n_k}: k \in \mathbb{N}\}) \geq \frac{\epsilon}{4}$. 
	\begin{corollary}\label{upperHausdorffcompatiblemetric}
		Let $X$ be a metrizable space and let $d, \rho$ be compatible metrics on $X$. Then the following are equivalent:
		\begin{enumerate}[(i)]
			\item $\tau_{H_d}^+ \supseteq \tau_{H_\rho}^+$ on $CL(X)$;
			\item for each $S \in \mathcal{P}_0(X)$ and $  f,g \in \mathcal{Z}^+$ such that $\inf_{x \in S}\{g(x)-f(x)\}>0$ and $B_\rho(S, g) \neq X$, then $B_\rho(S, f)$ is strict variationally $(\mathcal{P}_0(X)$-$d)$ included in $B_\rho(S, g)$;
			\item for each $S \in \mathcal{P}_0(X)$ and $ 0 < \alpha < \epsilon$ such that $B_\rho(S, \epsilon) \neq X$, then $B_\rho(S, \alpha)$ is strict variationally $(\mathcal{P}_0(X)$-$d)$ included in $B_\rho(S, \epsilon)$;
			\item the identity function $Id:(X,d) \to (X, \rho)$ is uniformly continuous.
		\end{enumerate}
	\end{corollary}
	\begin{proof}
The equivalence $(i) \Leftrightarrow (ii)$ follows from Theorem \ref{upperTsdcompatiblemetric} and $(ii) \Rightarrow (iii)$ is immediate.


$(iii) \Rightarrow (iv)$. Suppose $Id:(X,d) \to (X, \rho)$ is not uniformly continuous. Then we can find an $\epsilon > 0$ and sequences $(y_k)$ and $(z_k)$ such that $d(y_k,z_k) < \frac{1}{k}$ but $\rho(y_k,z_k) > 4\epsilon$. By the Efremovic lemma, there exist subsequences $(y_{k_n})$ and $(z_{k_n})$ such that $\rho(\{y_{k_n}: n \in \mathbb{N}\}, \{z_{k_n}: n \in \mathbb{N}\}) \geq \epsilon$. Take $A = \{y_{k_n}: n \in \mathbb{N}\}$. By $(ii)$, choose $D \in \mathcal{P}_0(X)$ and $f,g \in \mathcal{Z}^+$ such that $\inf_{x \in D}(g(x)-f(x)) = r>0$ and $$B_\rho(A, \frac{\epsilon}{4}) \subseteq B_d(D,f) \subseteq B_d(D, g) \subseteq B_\rho(A, \frac{\epsilon}{2}).$$ Then $B_d(A,r) \subseteq B_d(D, g)$ and since $d(y_{k}, z_{k}) < \frac{1}{k}$ $\forall k$, there exists $n_0 \in \mathbb{N}$ such that $\{z_{k_n}: n \geq n_0\} \subseteq B_d(D,g)$. But then $\epsilon \leq \rho(\{y_{k_n}: n \in \mathbb{N}\}, \{z_{k_n}: n \in \mathbb{N}\}) \leq \frac{\epsilon}{2}$. Hence, we arrive at a contradiction.

$(iv) \Rightarrow (ii)$. Take $S \in \mathcal{P}_0(X)$ and $f,g \in \mathcal{Z}^+$ with $\inf_{x \in S}(g(x)-f(x)) = \epsilon > 0$. Let $S' = B_\rho(S,f)$, then $S' \in \mathcal{P}_0(X)$. Choose $\delta > 0$ such that whenever $d(x,y) < \delta$, then $\rho(x,y) < \epsilon$. We show that $$B_\rho(S,f) \subseteq B_d(S', \frac{\delta}{3}) \subseteq B_d(S', \frac{\delta}{2})\subseteq B_\rho(S,g).$$  Clearly, $B_\rho(S,f) \subseteq B_d(S', \frac{\delta}{3})$. To see the last inclusion, take $y \in B_d(S', \frac{\delta}{2})$. Then $\exists$ $x' \in S'$ such that $d(x',y) < \frac{\delta}{2}$, so $\rho(x',y) < \epsilon$. Choose $x \in S$ such that $\rho(x',x) < f(x)$. Consequently, $\rho(x,y) < f(x) + \epsilon\leq g(x)$.    
	\end{proof}
	

\begin{lemma}\normalfont{(Proposition $4.3$, \cite{agarwal2024set})}\label{Bd(X)}
	Let $(X,d)$ be a metric space. If $S \in \mathcal{B}_d(X)$ and $f \in \mathcal{Z}^+$, then either $B_d(S,f)$ is equal to $X$, or lies in $\mathcal{B}_d(X)$. 
\end{lemma}	

	\begin{corollary}\label{AttocuhWetscompatiblemetric}
		Let $X$ be a metrizable space and let $d, \rho$ be two compatible metrics on $X$ that determine same bounded sets, that is, $\mathcal{B}_d(X)= \mathcal{B}_\rho(X)= \mathcal{B}(X)$ (say). Then the following are equivalent:
		\begin{enumerate}[(i)]
			\item $\tau_{AW_d}^+ \supseteq \tau_{AW_\rho}^+$;
			\item for each $B \in \mathcal{B}(X)$ and $f,g \in \mathcal{Z}^+$ such that $\inf_{x \in B}(g(x)-f(x))>0$ such that $B_\rho(B, g) \neq X$, then $B_\rho(B, f)$ is strict variationally $(\mathcal{B}(X)$-$d)$ included in $B_\rho(B, g)$;
			\item for each $B \in \mathcal{B}(X)$ and $ 0 < \alpha < \epsilon$ such that $B_\rho(B, \epsilon) \neq X$, then $B_\rho(B, \alpha)$ is strict variationally $(\mathcal{B}(X)$-$d)$ included in $B_\rho(B, \epsilon)$
			\item the Id: $(X,d) \to (X, \rho)$ is (strongly) uniformly continuous on $\mathcal{B}(X)$.
		\end{enumerate}
	\end{corollary}
	\begin{proof}
	The equivalence $(i) \Leftrightarrow (ii)$ follows from Theorem \ref{upperTsdcompatiblemetric}. The proof of $(ii) \Leftrightarrow (iii)$ is similar to the proof of $(ii) \Leftrightarrow (iii)$ of Corollary \ref{upperHausdorffcompatiblemetric}.
	
	$(iii) \Rightarrow (iv)$. Suppose  there exists $B \in \mathcal{B}(X)$ such that $Id: (X,d)\to (X, \rho) $ is not strongly uniformly continuous on $B$. So we can find an $\epsilon > 0$ and sequences $(y_k)$ and $(z_k)$ such that $\{y_k: k \in \mathbb{N}\}\subseteq B$ and $d(y_k, z_k)< \frac{1}{k}$ but $\rho(y_k, z_k) > 4\epsilon$. Proceeding now as in the proof of Corollary \ref{upperHausdorffcompatiblemetric}, we arrive at a contradiction.
	
	$(iv) \Rightarrow (ii)$. Take $B \in \mathcal{B}(X)$ and $f,g \in \mathcal{Z}^+$ with $\inf_{x \in B}(g(x)-f(x)) = \epsilon > 0$ such that  $B_{\rho}(B,g) \neq X$. Then by Lemma \ref{Bd(X)}, $B' = B_{\rho}(B,f) \in \mathcal{B}(X)$. By $(iv)$, choose $\delta > 0$ such that whenever $\{x',y\}\cap S \neq \emptyset$ and $d(x',y) < \delta$, then $\rho(x',y) < \epsilon$. The rest of the proof is similar to the proof of $(iv) \Rightarrow (ii)$ of Corollary \ref{upperHausdorffcompatiblemetric}. 
	\end{proof}
	
	The following natural question arises from the relation between strong uniform continuity of $Id:(X,d) \to (X, \rho)$ on $\mathcal{B}_d(X)$ (resp. $\mathcal{P}_0(X)$) and the inclusion $\tau_{{AW}_d}^+ \supseteq \tau_{AW_\rho}^+$ (resp. $\tau_{H_{d}}^+ \supseteq \tau_{H_\rho}^+$) in Corollary \ref{AttocuhWetscompatiblemetric} (resp. Corollary \ref{upperHausdorffcompatiblemetric}): whether a similar relation holds for an arbitrary bornology $\mathcal{S}$? The answer is partially affirmative. Interestingly, in the case of $\tau_{\mathcal{S},d}^+$-convergence, the strong uniform continuity of the identity function on $\mathcal{S}$ as well as on $\mathcal{S}_{\mathcal{Z}^+}^d$ play a central role.  Recall that $\mathcal{S}_{\mathcal{Z}^+}^d$ denotes the set of all proper functional enlargements of members of $\mathcal{S}$.

 \begin{proposition}\label{strongunifromcontinuitysuffcondition}
 	Let $X$ be a metrizable space and let $\mathcal{S}$ be a bornology on $X$. Suppose $d$ and $\rho$ are two compatible metrics on $X$. Consider the following statements:
 	\begin{enumerate}[(i)]
 		\item for $S \in \mathcal{S}$ and  $f,g \in \mathcal{Z}^+$ whenever $\inf_{x \in S}(g(x)-f(x))>0$, then $B_\rho(S,f)$ is strict variationally $(\mathcal{S}$-$d)$ included in $B_\rho(S,g)$;
 		\item the identity function $Id: (X,d) \to (X, \rho)$ is strongly uniformly continuous on $\mathcal{S}$. 
 		 \end{enumerate}
 	Then $(i) \Rightarrow (ii)$ holds.
 	  \end{proposition}
 	  \begin{proof}
 	  	Suppose $Id:(X,d) \to (X, \rho)$ is not strongly uniformly continuous on $\mathcal{S}$. Then there exists $S_0 \in \mathcal{S}$ such that $Id:(X,d) \to (X, \rho)$ is not strongly uniformly continuous on $S_0$. So we can find an $\epsilon >0$ and sequences $(x_n) \subseteq S_0$ and $(y_n)$ in $X$ such that $d(x_n,y_n) < \frac{1}{n}$ but $\rho(x_n,y_n) \geq 4\epsilon$ $\forall$ $n$. Imitating the proof of Corollary \ref{upperHausdorffcompatiblemetric}, we arrive at a contradiction.  
 	  \end{proof}

 	  \begin{theorem}\label{bornologicalconvergncesuff}
 	  	Let $X$ be a metrizable space and let $\mathcal{S}$ be a bornology on $X$. Suppose $d, \rho$ are two compatible metrics on $X$. Consider the following statements:
 	  	\begin{enumerate}[(i)]
 	  		\item $\tau_{\mathcal{S},d}^+ \supseteq \tau_{\mathcal{S},\rho}^+$ on $CL(X)$;
 	  		\item $(\mathcal{S}-d) \geqslant (\mathcal{S}-\rho)$ on $CL(X)$;
 	  		\item $(\mathcal{S}$-$d)^+ \geqslant (\mathcal{S}$-$\rho)^+$ on $CL(X)$. 
 	  	\end{enumerate}
 	  	Then $(i) \Rightarrow (ii) \Leftrightarrow (iii)$ holds.
 	  \end{theorem}
 	  \begin{proof}$(i)\Rightarrow (ii)$.
 	  	By Theorem \ref{upperTsdcompatiblemetric}, $(i)$ is equivalent to say that for each $S \in \mathcal{S}$ and $f,g \in \mathcal{Z}^+$ such that $\inf_{x \in S}(g(x)-f(x))>0$, whenever $B_\rho(S,g) \neq X$, then $B_\rho(S,f)$ is strict variationally $(\mathcal{S}-d)$ included in $B_\rho(S,g)$. By Proposition \ref{strongunifromcontinuitysuffcondition}, the identity function $Id: (X,d) \to (X, \rho)$ is strongly uniformly continuous on $\mathcal{S}$. Consequently, invoking Theorem $3.13$ of \cite{PbciAWc}, we have $(\mathcal{S}$-$d) \geqslant (\mathcal{S}$-$\rho)$ on $CL(X)$.
 	  	
 	  	$(ii) \Rightarrow (iii)$. Let $(A_{\lambda})$ be a net in $CL(X)$ which is $(\mathcal{S}-d)^+$-convergent to $A$. Then it is easy to verify that $(A_{\lambda} \cup A)$ is $(\mathcal{S}-d)$-convergent to $A$. By $(ii)$, $(A_{\lambda} \cup A)$ is $(\mathcal{S}-\rho)$-convergent to $A$. So $A_{\lambda}\cup A \xrightarrow{(\mathcal{S}-\rho)^+} A$. Consequently, it is easy to see that $A_{\lambda}\xrightarrow{(\mathcal{S}-\rho)^+}A$ on $CL(X)$.
 	  	
 	  	$(iii) \Rightarrow (ii)$. By Theorem $3.13$ of \cite{PbciAWc}, $(\mathcal{S}$-$d)^+ \geqslant (\mathcal{S}$-$\rho)^+$ is equivalent to the fact that identity function $Id:(X,d) \to (X, \rho)$ is strongly uniformly continuous on $\mathcal{S}$. Further, invoking Theorem $3.1$ of \cite{Suc}, it is enough to show that for each $S \in \mathcal{S}$ and $A \in CL(X)$ whenever $D_d(S,A) =0$, then $D_{\rho}(S,A) =0$.
 	  	To see this, let $S \in \mathcal{S}$ and $A \in CL(X)$ such that $D_d(S,A) = 0$ but $D_{\rho}(S,A) = r > 0$. Then $S \cap B_d(A, \frac{1}{n}) = \emptyset$. Define $A_n = \{a_n\}$ for each $n \in \mathbb{N}$. So $(A_n)$ is a sequence in $CL(X)$.  It is easy to see that $(A_n)$ is $(\mathcal{S}$-$d)^+$-convergent to $A$. On the hand, observe that $S \cap B_d(A,r) = \emptyset$. Consequently, $(A_n)$ fails $(\mathcal{S}$-$\rho)^+$-convergent to $A$. 
 	  \end{proof}

 	  The implication $(iii) \Rightarrow (i)$ in Theorem \ref{bornologicalconvergncesuff} need not hold in general. We show it by the following example.
 	  \begin{example}
Let $X = \mathbb{R}$ and let $d = d_{0,1}$, the $0$-$1$ discrete metric on $\mathbb{R}$, and $\rho = d_u$, the usual metric on $\mathbb{R}$. Clearly, the Identity function $Id: (\mathbb{R}, d_{0,1}) \to (\mathbb{R},d_u)$ is continuous. So by Proposition $1.2$ of \cite{Suc}, $Id$ is strongly uniformly continuous on $\mathcal{F}(X)$. Invoking Theorem $3.13$ of \cite{PbciAWc}, we get $(\mathcal{F}(X)-d_{0,1})^+ \geqslant (\mathcal{F}(X)-d_u)^+$. We now show that $\tau_{W_{d_{0,1}}}^+\ngeqslant \tau_{W_{d_u}}^+$. Suppose $A_n = \{\frac{1}{n}\}$ for each $n \in \mathbb{N}$ and $A = \{1\}$. Then one can easily verify that $(A_n)$ is $\tau_{W_{d_{0,1}}}^+$-convergent to $A$ but is not $\tau_{W_{d_u}}^+$-convergent to $A$. \qed      	  	
 	  \end{example}

 	  The following result gives a sufficient condition under which all the statements of Theorem \ref{bornologicalconvergncesuff} become equivalent. 
 	  
 	  \begin{proposition}\label{Tsdbornological}
 	  	Let $X$ be a metrizable space and let $\mathcal{S}$ be a bornology on $X$. Suppose $d$ and $\rho$ are two compatible metrics on $X$ are such that $\mathcal{S} = \mathcal{S}^d_{\mathcal{Z}^+} =\mathcal{S}^\rho_{\mathcal{Z}^+}$. Then $\tau_{\mathcal{S},d}^+ \supseteq \tau_{\mathcal{S},\rho}^+$ if and only if $(\mathcal{S}-d) \geqslant (\mathcal{S}-\rho)$ on $CL(X)$.
 	  \end{proposition}
 	  \begin{proof}
 	  	It is enough to show that $(\mathcal{S}-d) \geqslant (\mathcal{S}-\rho)$ implies $\tau_{\mathcal{S},d}^+ \supseteq \tau_{\mathcal{S},\rho}^+$. Let $S\in \mathcal{S}$ and $f,g \in \mathcal{Z}^+$ such that $\inf_{x \in S}(g(x)-f(x))  =\epsilon>0$
 	  	and $B_\rho(S,g) \neq X$. Since $\mathcal{S} = \mathcal{S}^d_{\mathcal{Z}^+} =\mathcal{S}^\rho_{\mathcal{Z}^+}$, take $S' = B_\rho(S,f) \in \mathcal{S}$. Since $(\mathcal{S}-d) \geqslant (\mathcal{S}-\rho)$, by Theorem $3.13$ of \cite{PbciAWc}, the identity function $Id: (X,d) \to (X, \rho)$ is strongly uniformly continuous on $\mathcal{S}$. Choose $\delta > 0$ such that whenever $x' \in S'$ and $y \in X$ such that $d(x', y) < \delta$, then $\rho(x',y) < \epsilon$. So $$B_\rho(S,f)\subseteq B_d(S', \frac{\delta}{3})\subseteq B_d(S', \frac{\delta}{2}) \subseteq B_\rho(S, g).$$ Thus, by Theorem \ref{upperTsdcompatiblemetric}, $\tau_{\mathcal{S},d}^+ \supseteq \tau_{\mathcal{S},\rho}^+$ on $CL(X)$.   
 	  \end{proof}
 	  
 	  Alternatively, Proposition \ref{Tsdbornological} also follows from  Theorem $6.4$ of \cite{agarwal2024set} as when $\mathcal{S} = \mathcal{S}^d_{\mathcal{Z}^+} =\mathcal{S}^\rho_{\mathcal{Z}^+}$, we have $\tau_{\mathcal{S},d}^+ = (\mathcal{S}-d)^+$ and $\tau_{\mathcal{S},\rho}^+ = (\mathcal{S}-\rho)^+$.  
 
 Now we establish a sufficient condition for the equivalence $\tau_{\mathcal{S},d}^+ = \tau_{\mathcal{S},\rho}^+$ in terms of strong uniform continuity. Let $X$ be a metrizable space and let $\mathcal{S}$ be a bornology on $X$. For two compatible metrics $d,\rho$ on $X$, we say $\mathcal{S}^d_{\mathcal{Z}^+}\geqq \mathcal{S}^\rho_{\mathcal{Z}^+}$ if for each $S \in \mathcal{S}$, $f \in \mathcal{Z}^+$, and $\epsilon > 0$, $\exists S' \in \mathcal{S}$ and $h \in \mathcal{Z}^+$ such that $B_\rho(S,f) \subseteq B_d(S',h) \subseteq B_\rho(S,f+\epsilon)$. For example, for any two equivalent norms $\| \cdot \|_1$ and $\| \cdot \|_2$ on a vector space $X$, we have $\mathcal{S}^{\| \cdot \|_1}_{\mathcal{Z}^+} \geqq \mathcal{S}^{\| \cdot \|_2}_{\mathcal{Z}^+}$, where $\mathcal{S}$ is the collection of all nonempty bounded subsets of $X$. 
 
 \begin{proposition}\label{Tsdcomatbilmetricsuff}
 	Let $X$ be a metrizable space and let $\mathcal{S}$ be a bornology on $X$.  Suppose $d$ and $\rho$ are two compatible metrics on $X$ such that $(X,d)$ is an a.c.m.s. and $\mathcal{S}^d_{\mathcal{Z}^+} \geqq \mathcal{S}^\rho_{\mathcal{Z}^+} $. Consider the following statements:
 	\begin{enumerate}[(i)]
 		\item the identity function, $Id:(X,d) \to (X,\rho)$ is strongly uniformly continuous on $\mathcal{S}^d_{\mathcal{Z}^+}$;
 		\item  for $S \in \mathcal{S}$ and  $f,g \in \mathcal{Z}^+$ whenever $\inf_{x \in S}(g(x)-f(x))>0$ and $B_{\rho}(S,g) \neq X$, then $B_\rho(S,f)$ is strict variationally $(\mathcal{S}$-$d)$ included in $B_\rho(S,g)$.   
 	\end{enumerate} 	
 	Then $(i) \Rightarrow (ii)$ holds.
 \end{proposition}
 \begin{proof}
 	Let $S \in \mathcal{S}$, $f,g \in \mathcal{Z}^+$ with $\inf_{x \in S}(g(x)-f(x)) = \epsilon > 0$ and $B_{\rho}(S,g) \neq X$. Choose $S' \in \mathcal{S}$ and $h \in \mathcal{Z}^+$
 	such that $B_\rho(S,f) \subseteq B_d(S',h) \subseteq B_\rho(S, f+\frac{\epsilon}{2})$. By $(i)$, $\exists$ $\delta > 0$ such that whenever $d(x,y) < \delta$ and $x \in B_d(S',h)$, then $\rho(x,y) < \frac{\epsilon}{2}$. We now show that $$B_\rho (S,f) \subseteq B_d(S',h) \subseteq B_d(S', h+\frac{\delta}{2}) \subseteq B_\rho(S,g).$$	Take $y \in B_d(S', h+\frac{\delta}{2})$. Since $(X,d)$ is an a.c.m.s., $\exists$ $x' \in B_d(S',h)$  such that $d(x',y) < \frac{\delta}{2}$. Then $\rho(x',y) < \frac{\epsilon}{2}$. Also $B_d(S',h) \subseteq B_\rho(S, f+\frac{\epsilon}{2})$ so $y \in B_\rho(S,g)$. Thus, $B_d(S',h+\frac{\delta}{2}) \subseteq B_\rho(S,g)$. \end{proof}
 
 The following result follows from Theorem \ref{upperTsdcompatiblemetric} and Proposition \ref{Tsdcomatbilmetricsuff}.
 \begin{theorem}
 	Let $X$ be a metrizable space and let $\mathcal{S}$ be a bornology on $X$. Suppose $d$ and $\rho$ are two compatible metrics on $X$ such that $(X,d)$ and $(X, \rho)$ are a.c.m.s. and $\mathcal{S}^d_{\mathcal{Z}^+}\geqq \mathcal{S}^\rho_{\mathcal{Z}^+}$, $\mathcal{S}^\rho_{\mathcal{Z}^+}\geqq \mathcal{S}^d_{\mathcal{Z}^+}$ hold. Then $\tau_{\mathcal{S},d}^+ = \tau_{\mathcal{S},\rho}^+$ if the identity functions $Id: (X,d) \to (X, \rho)$ and  $Id: (X, \rho) \to (X,d)$ are strongly uniformly continuous on $\mathcal{S}^d_{\mathcal{Z}^+}$ and on $\mathcal{S}^\rho_{\mathcal{Z}^+}$, respectively.  
 \end{theorem}
 	
 	\section{Infima of hyperspace convergences}
 	In the final section of this paper, we examine the problem of infima of hyperspace convergences corresponding to a given family of equivalent metrics on a metrizable space $X$. Particular focus is placed  on the classical Attouch-Wets and Wijsman convergences. We show that a metrizable space $X$ admits a compatible Heine-Borel metric if and only if the collection $\{\tau_{AW_\rho}^+: \rho \text{ is a compatible metric on } X\}$ has minima. Furthermore, for a metric space $(X,d)$ and any two bornologies $\mathcal{S}$ and $\mathcal{B}$, we examine the inclusion $\tau_{\mathcal{S},d}^+ \subseteq \tau_{\mathcal{B},\rho}^+$ for any  uniformly equivalent metric $\rho$. Consequently, it is proved that a metrizable space $X$ is separable if and only if $X$ admits a compatible metric $d$ such that the collection $\{\tau_{W_\rho}^+: \rho \text{ is uniformly equivalent to }d\}$ has $\tau_{AW_{d}}^+$ as infimum.
 	

We begin by examining the conditions on $\mathcal{S}$ and on the metric $d$ for which $\tau_{\mathcal{S},d}^+$ is the minimum of the collection $\{\tau_{\mathcal{S},\rho}^+: \rho \text{ is uniformly equivalent to } d \}$.

 	The following lemma is a reformulation of Theorem \ref{Miss theorem for upper Tsd} for the metric $\rho = \min\{d,M\}$ on $(X,d)$, where $M> 0$. The proof technique being similar to that of Theorem \ref{Miss theorem for upper Tsd}, we omit its proof.
 	\begin{lemma}\label{Misstyperho}
 		Let $(X,d)$ be a metric space and let $\mathcal{S}$ be a bornology on $X$. Suppose $\rho = \min\{M,d\}$ for some $M> 0$, $(A_{\lambda})$ is a net in $CL(X)$ and $A \in CL(X)$. Then the following statements are equivalent:
 		\begin{enumerate}[(i)]
 			\item $(A_{\lambda})$ is $\tau_{\mathcal{S},\rho}^+$-convergent to $A$;
 			\item for $S \in \mathcal{S}$ and $f,g \in \mathcal{Z}^+$ such that $\inf_{x \in S}(g(x)-f(x))>0$ and $g(x) \leq M$ $\forall x \in S$ whenever $A$ misses $B_d(S,g)$, then eventually $(A_{\lambda})$ misses $B_d(S,f)$.
 		\end{enumerate}
 	\end{lemma}
 	\begin{theorem}\label{Tsdunirformlyequimetric}
 		Let $(X,d)$ be a metric space and let $\mathcal{S}$ be a bornology on $X$. Consider the following statements:
 		\begin{enumerate}[(i)]
 			\item $\tau_{\mathcal{S},d}^+ \subseteq \tau_{\mathcal{S},\rho}^+$ on $CL(X)$ for each metric $\rho$ that is uniformly equivalent to $d$;
 			\item for $S \in \mathcal{S}$ and $f,g \in \mathcal{Z}^+$ such that $\inf_{x \in S}(g(x)-f(x))>0$ and $B_d(S,g) \neq X$, we have $B_d(S,f) \in \mathcal{S}_d^*$, that is $B_d(S,f)$ is weakly $\mathcal{S}$-totally bounded with respect to $d$. 
 		\end{enumerate}
 		Then $(i) \Rightarrow (ii)$ holds. If in $(ii)$, we have $B_d(S,f) \in \mathcal{S}^d_*$, then the converse also holds.
 	\end{theorem}
 	\begin{proof}
 		$(i)\Rightarrow (ii)$. Suppose $(ii)$ fails. Then we can find an $S_{0} \in \mathcal{S}$ and $f,g \in \mathcal{Z}^+$ with $\inf_{x \in S}(g(x)-f(x))>0$ satisfying $B_d(S_0,g) \neq X$ and $B_d(S_0,f)$ is not weakly $\mathcal{S}$-totally bounded. So we can find an $r >0$ such that for each $S \in \mathcal{S}$, $B_d(S_0,f) \nsubseteq B_d(S,2r)$. For each $S \in \mathcal{S}$, choose $x_S \in B_d(S_0,f) \setminus B_d(S,2r)$. Direct $\mathcal{S}$ by the set inclusion $\subseteq$. Define $A = X\setminus B_d(S_0,g)$ and for every $S \in \mathcal{S}$, $A_S = A\cup\{x_S\}$. Then $(A_S)$ is a net in $CL(X)$ and $A\in CL(X)$. Let $\rho = \min\{r,d\}$. Then $\rho$ is an uniformly equivalent metric to $d$ on $X$. We show that $(A_S)$ is $\tau_{\mathcal{S},\rho}^+$-convergent to $A$ but it is not $\tau_{\mathcal{S},d}^+$-convergent to $A$.  
 		
 		To see $\tau_{\mathcal{S},\rho}^+$-convergence, we use Lemma \ref{Misstyperho}. Let $S_0 \in \mathcal{S}$ and $p,q \in \mathcal{Z}^+$ such that $\inf_{x \in S_0}(q(x)-p(x))>0$, $q(x) \leq r$ $\forall x \in S_0$, and $A \cap B_d(S_0,q) = \emptyset$. We show that $(A_S)$ eventually misses $B_d(S_0,p)$.
 		If $S \supseteq S_0$, then $2r \leq d(x_{S},S) \leq d(x_{S}, S_0)$. Consequently, $A_S\cap B_d(S_0,p) = \emptyset$. 
 		
 		Next, we show $(A_S)$ is not $\tau_{\mathcal{S},d}^+$-convergent to $A$. Observe that $A \cap B_d(S_0,g) = \emptyset$ but for every $S \in \mathcal{S}$, we have $A_{S} \cap B_d(S_0,f) \neq \emptyset$. Thus, by Theorem \ref{Miss theorem for upper Tsd}, $(A_S)$ is not $\tau_{\mathcal{S},d}^+$-convergent to $A$.
 		
 		Conversely, suppose in $(ii)$, $B_{d}(S,f)$ is $\mathcal{S}$-totally bounded with respect to $d$. Let $\rho$ be any uniformly equivalent metric to $d$.  Suppose $(A_{\lambda})$ is a net which is not $\tau_{\mathcal{S},d}^+$-convergent to $A$ in $CL(X)$. We show that $(A_{\lambda})$ is not $\tau_{\mathcal{S},\rho}^+$-convergent to $A$. Take $S_0 \in \mathcal{S}$ and $f,g \in \mathcal{Z}^+$ such that $\inf_{x \in S_0}(g(x)-f(x)) =r>0$ and $A \cap B_d(S_0,g) = \emptyset$ but $A_{\lambda} \cap B_d(S_0,f) \neq \emptyset$ frequently. Since $d$ and $\rho$ are uniformly equivalent, $\exists$ $\epsilon >0$ such that whenever $\rho(x,y)<\epsilon$, then $d(x,y)<r$ and $\exists$ $\delta > 0$ such that whenever $d(x,y) < \delta$, then $\rho(x,y) < \frac{\epsilon}{2}$.  Given $B_d(S_0,f)$ is $\mathcal{S}$-totally bounded, $\exists$ $S \in \mathcal{S}$ with $S \subseteq B_d(S_0,f)$ satisfying $B_d(S_0,f) \subseteq B_d(S, \delta) \subseteq B_\rho(S, \frac{\epsilon}{2})$. So $A_{\lambda} \cap B_{\rho}(S, \frac{\epsilon}{2}) \neq \emptyset$ frequently. Also $B_\rho (S, \epsilon) \subseteq B_d(S,r)$. Consequently, $B_\rho(S,\epsilon) \subseteq B_d(S_0,g)$, so  $A \cap B_\rho(S, \epsilon) = \emptyset$. Thus, $(A_{\lambda})$ is not $\tau_{\mathcal{S},\rho}^+$-convergent to $A$.           
 	\end{proof} 
 	Note that by Theorem $13.13$ of \cite{beer2023bornologies} whenever $\mathcal{S}_*^d$ forms a bornology on $X$, that is, when $\mathcal{S}_*^d  =\mathcal{S}_d^*$, then statements $(i)$ and $(ii)$ in Theorem \ref{Tsdunirformlyequimetric} become equivalent.
 	
 	Under the assumption of almost convexity on the base space $(X,d)$, we get the following result. 
 	\begin{corollary}
 		Let $(X,d)$ be an a.c.m.s. and let $\mathcal{S}$ be a bornology on $X$. Consider the following statements:
 		\begin{enumerate}[(i)]
 			\item $\tau_{\mathcal{S},d}^+ \subseteq \tau_{\mathcal{S},\rho}^+$ on $CL(X)$ for each uniformly equivalent metric $\rho$ on $(X,d)$;
 			\item each non-dense member of $\mathcal{S}_{\mathcal{Z}^+}^d$ is weakly $\mathcal{S}$-totally bounded with respect to $d$. 		\end{enumerate}
 		Then $(i) \Rightarrow (ii)$ holds. If in $(ii)$, each non-dense member of $\mathcal{S}_{\mathcal{Z}^+}^d$ is in $\mathcal{S}_*^d$, then the converse also holds.
 	\end{corollary}

 	\begin{corollary}
 		Let $(X,d)$ be a metric space. Then the following statements are equivalent:
 		\begin{enumerate}[(i)]
 			\item $\tau_{W_d} \subseteq \tau_{W_\rho}$ for any uniformly equivalent metric $\rho$ on $X$;
 			\item each proper $d$-closed ball is $d$-totally bounded.  
 		\end{enumerate}
 	\end{corollary}
 	
 	The above corollary can also be obtained by applying Lemma $4.1$ and Theorem $4.3$ of \cite{LechickiWijsman}.
 	
 	By slightly modifying the proof of Theorem \ref{upperTsdcompatiblemetric}, one can prove the following result. 	
 	\begin{theorem}\label{Tsddiffbornologiesmetrics}
 		Let $X$ be a metrizable space  and let $d, \rho$ be two compatible metrics on $X$. Suppose $\mathcal{S}$ and $\mathcal{C}$ are two bornologies on $X$. Then the following statements are equivalent:
 		\begin{enumerate}[(i)]
 			\item  $\tau_{\mathcal{S},d}^+ \supseteq \tau_{\mathcal{C}, \rho}^+$ on $CL(X)$;
 			\item  for each $C \in \mathcal{C}$ and $f,g \in \mathcal{Z}^+$ whenever $\inf_{x \in C}(g(x)-f(x))>0$ and $B_{\rho}(C,g) \neq X$, then $B_{\rho}(C,f)$ is strict variationally $(\mathcal{S}$-$d)$ included in $B_{\rho}(C,g)$;
 			\item for each pair $(A,C)$ of nonempty proper subsets of $X$ such that whenever $A$ is strict variationally $(\mathcal{C}$-$\rho)$ included in $C$, we have $A$ is strict variationally $(\mathcal{S}$-$d)$ included in $C$.
 			\end{enumerate} 	\end{theorem}
 		
In particular, for an arbitrary bornology $\mathcal{S}$, the previous theorem implies that $\tau_{W_d}^+ \supseteq \tau_{\mathcal{S}, \rho}^+$ if and only if whenever $A$ is strict variationally ($\mathcal{S}$-$\rho$) included in $C$, then $A$ is strictly $d$-included in $C$. Thus, if $\rho =d$, then the topology $\tau_{\mathcal{S},d}^+$ has a weak representation by the family $\{d(x, \cdot): x \in X\}$ of distance functionals provided the notions of strictly $d$-inclusion and strict variationally  ($\mathcal{S}$-$d)$ inclusion are equivalent for proper nonempty subsets of $X$. The converse also holds in this case.
 		\begin{corollary}$($Theorem 5.2, \cite{beer2013gap}$)$
 		Let $(X,d)$ be a metric space. Then the following statements are equivalent:
 		\begin{enumerate}[(i)]
 			\item $\tau_{{AW}_d}^+ = \tau_{H_{d}}^+$;
 			\item $(X,d)$ is bounded.
 		\end{enumerate} 
 	\end{corollary}
 	\begin{proof}
 		$(i) \Rightarrow (ii)$.	Let $p \in X$ and $A = X \setminus B_d(p, \epsilon)$ for some $\epsilon > 0$. If $A = \emptyset$, then we are done. Otherwise, we show that $A$ is bounded. Since $\tau_{AW_d}^+ = \tau_{H_{d}}^+$, by Theorem \ref{Tsddiffbornologiesmetrics}, we can find $S \in \mathcal{B}_d(X)$ and $f,g \in \mathcal{Z}^+$ with $\inf_{x \in S}(g(x) - f(x)) > 0$ satisfying $$ B_d(A, \frac{\epsilon}{4}) \subseteq B_d(S,f) \subseteq B_d(S,g) \subseteq B_d(A, \frac{\epsilon}{2}).$$ But by Lemma \ref{Bd(X)}, $B_d(S,g) \in \mathcal{B}_d(X)$. So $A$ is bounded. Consequently, $(X,d)$ is bounded. 
 		
 		$(ii)\Rightarrow (i)$. It is  immediate from Theorem \ref{Tsddiffbornologiesmetrics}. 
 	\end{proof}
 
  It is well-known that a metrizable space is compact if and only if the  Wijsman topologies corresponding to all compatible metrics coincide (Corollary $5.6$, \cite{beer1992suprema}). Also a metric space $(X,d)$ has nice closed balls (proper closed balls are compact) if and only if the family of all Wijsman topologies over all equivalent metrics to $d$ has $\tau_{W_d}$ as its minimum element \cite{infinmahyperspace}. The next result answers a similar problem for the upper Attouch-Wets convergence. 
 		Recall that a compatible metric $d$ on $X$ is called a \textit{Heine-Borel metric} if all closed and bounded subsets of $(X,d)$ are compact. 
 		
 		\begin{theorem}\label{AWequivalentminimum}
 			Let $X$ be a metrizable space. Then the following statement are equivalent: 
 			\begin{enumerate}[(i)]
 				\item $X$ admits a  compatible metric $d$ such that $$\tau_{{AW}_d}^+ = \min \{\tau_{AW_\rho}^+: \rho \text{ is a compatible metric on } X\};$$
 				\item $X$ admits a compatible Heine-Borel metric.
 			\end{enumerate}
 		\end{theorem}
 		\begin{proof}
 			$(i) \Rightarrow (ii)$. Suppose $X$ admits a compatible metric $d$ such that $\tau_{{AW}_d}^+ = \min \{\tau_{AW_\rho}^+: \rho \text{ is a compatible metric on } X\}$.  Let $B$ be a proper closed and bounded set in $(X,d)$. Take $p \in X \setminus B$ and $\epsilon > 0$ such that $d(p, B) > 3\epsilon$. By $(i)$ and Theorem \ref{Tsddiffbornologiesmetrics}, $B_d(B, \epsilon)$ is strict variationally $(\mathcal{B}_\rho(X)$-$\rho)$ included in $B_d(B, 2\epsilon)$ for every compatible metric $\rho$ on $X$. Therefore, $B \in \mathcal{B}_\rho(X)$ for each compatible metric $\rho$ on $X$. Consequently, by Theorem $6.8$ of \cite{beer2023bornologies}, $B$ is compact. Thus, each proper closed and bounded subset of $(X,d)$ is compact. If $(X,d)$ is unbounded, then $d$ is the required Heine-Borel metric. Otherwise, fix $x_0 \in X$. Choose $0< r < 1$ such that $B_d(x_0, 2r) \neq X$. Then $X = \overline{B_d(x_0, r)} \cup (X\setminus B_d(x_0,r))$ and both $\overline{B_d(x_0, r)}$ and $X\setminus B_d(x_0,r)$ are compact, thus $X$ is compact.

 			$(ii) \Rightarrow (i)$.  Suppose $X$ admits a compatible Heine-Borel metric, say $d$. Then $\mathcal{B}_d(X)  = \mathcal{K}(X)$. By Theorem \ref{Tsddiffbornologiesmetrics}, it is enough to show that for every compatible metric $\rho$ on $X$, whenever $S \in \mathcal{B}_d(X)$ and $f,g \in \mathcal{Z}^+$ such that $\inf_{x \in S}(g(x)-f(x))>0$ and $B_d(S,g) \neq X$, then $B_d(S,f)$ is strict variationally $(\mathcal{B}_\rho(X)$-$\rho)$ included in $ B_d(S,g)$.  By Lemma \ref{Bd(X)}, $B_d(S,g) \in \mathcal{B}_d(X)$. Since $\overline{B_d(S,f)} \subseteq B_d(S,g)$,  $S' = \overline{B_d(S,f)}$ is compact, thus $S' \in \mathcal{B}_\rho(X)$ for every compatible metric $\rho$ on $X$. Clearly, $B_d(S,g)$ is open in $(X, \rho)$ for every compatible metric $\rho$ on $X$, so for some $\epsilon > 0$ we have $$B_d(S,f) \subseteq B_\rho(S', \epsilon) \subseteq B_\rho(S',2\epsilon) \subseteq B_d(S,g).$$     
 		\end{proof}
 The above theorem can also be deduced using Corollary 2.5 in \cite{beer2013gap}. However, in \cite{beer2013gap}, the authors' approach for the Attouch-Wets convergence is through gap functionals.
 		

 	\begin{remark}\begin{enumerate}[(i)]
 			\item It is well-known that there exists no compatible Heine-Borel metric on an infinite dimensional normed linear space. Using this fact and Theorem \ref{AWequivalentminimum}, it can be deduced that on a normed linear space $X$, the collection of $\{\tau_{AW_{\rho}}^+: \rho \text{ is a compatible metric on }X\}$ has minimum element if and only if $X$ is finite dimensional.
 			
 			\item From Theorem 12 of \cite{costantini1993metrics} and Theorem \ref{AWequivalentminimum}, it also follows that for a metrizable space $X$, whenever there exists a compatible metric $d$ on $X$ such that $\tau_{{AW}_d}^+$ is the minimum of all $\tau_{AW_\rho}^+$s' such that $\rho$ is equivalent to $d$, then $\tau_{W_d}^+$ is the minimum of all $\tau_{W_\rho}^+$s' such that $\rho$ is equivalent to $d$.

 		\end{enumerate} 
 		\end{remark}
 		We now give a consequence of Theorem \ref{AWequivalentminimum} related to function spaces. Recall that by identifying a continuous function $f$ between two metric spaces with its graph, we can consider $f$ as a closed subset of the product space. In the next result, $d_u$ denotes the usual metric on $\mathbb{R}$.
 	\begin{corollary}
 		 Let $X$ be a metrizable space. Suppose the collection $\{\tau_{AW_{\rho}}^+: \rho \text{ is a compatible metric on }X\}$ has a minimum element, then there exists a compatible metric $d$ on $X$ such that each sequence of continuous functions from $X$ to $[0,1]$ is $\tau_{{AW}_{d\times d_u}}$-convergent to a continuous function $f:X \to \mathbb{R}$,  if and only if it is uniformly convergent to $f$ on $\mathcal{B}_d(X)$.\end{corollary}
 	 \begin{proof}
 		Suppose the collection $\{\tau_{AW_{\rho}}^+: \rho \text{ is a compatible metric on }X\}$ has minimum on $CL(X)$. Then by Theorem \ref{AWequivalentminimum}, there exists a compatible metric $d$ such that all closed and bounded subsets of $(X,d)$ are compact. We show that $d$ is the required metric on $X$. Let $(f_n)$ be a sequence of continuous functions from $X$ to $[0,1]$ and $f$ be a real-valued continuous function on $X$ such that $f_n \xrightarrow{\tau_{{AW}_{d\times d_u}}} f$. We show that $f_n \to f$ uniformly on $\mathcal{B}_d(X)$. Let $B \in \mathcal{B}_d(X)$ such that $B$ is closed. Then $B$ is compact. So $f$ is uniformly continuous on $B$ and clearly $\{f_n: n \in \mathbb{N}\}$ is uniformly bounded on $B$. Consequently, by Theorem 4.1  of \cite{beer1991awucb}, $(f_n)$ is uniformly convergent to $f$ on $B$.    
 	\end{proof}
 		In order to study relation between $\tau_{{AW}_d}^+$ and $\tau_{W_\rho}^+$ with respect to two compatible metrics $d, \rho$ on a metrizable space $X$, we first prove a general result.
 		\begin{proposition}\label{AWTsdcomparison}
 			Let $X$ be a metrizable space and let $\mathcal{S}$ be a bornology on $X$. Suppose $d,\rho$ are two compatible metrics on $X$ such that $(X, \rho)$ is an a.c.m.s. Consider the following statements:
 			\begin{enumerate}[(i)]
 				\item $\tau_{{AW}_d}^+ \supseteq \tau_{\mathcal{S},\rho}^+$ on $CL(X)$;
 				\item for $S \in \mathcal{S}$, $f \in \mathcal{Z}^+$ either  $B_{\rho}(S,f)$ is dense in $X$ or $B_\rho(S,f) \in \mathcal{B}_d(X)$.  
 			\end{enumerate}
 			Then $(i) \Rightarrow (ii)$ holds. The converse holds if in $(ii)$, the identity function $Id:(X,d) \to (X, \rho)$ is strongly uniformly continuous on $\mathcal{S}_{\mathcal{Z}^+}^\rho$. 
 		\end{proposition}
\begin{proof}$(i)\Rightarrow (ii)$.
	Let $S \in \mathcal{S}$ and $f \in \mathcal{Z}^+$ such that $\overline{B_{\rho}(S,f)} \neq X$. Then $\exists$ $p \in X\setminus \overline{B_{\rho}(S,f)}$. Choose $\epsilon >0$ such that $\rho(p, \overline{B_{\rho}(S,f)}) > 3\epsilon$. Take $g = f + \epsilon$, then $g \in \mathcal{Z}^+$. Since $(X, \rho)$ is an a.c.m.s., $B_\rho(S,g) \neq X$. By $(i)$ and Theorem \ref{Tsddiffbornologiesmetrics}, $B_\rho(S,f)$ is strict variationally $(\mathcal{B}_d(X)$-$d)$ included in $B_\rho(S,g)$. Consequently, by Lemma \ref{Bd(X)}, $B_\rho(S,f) \in \mathcal{B}_d(X)$.

Conversely, suppose $(ii)$ holds and in addition $Id:(X,d) \to (X, \rho)$ is strongly uniformly continuous on $\mathcal{S}_{\mathcal{Z}^+}^\rho$. We show that $(i)$ holds. Let $S \in \mathcal{S}$ and $f,g \in \mathcal{Z}^+$ such that $\inf_{x \in S}(g(x)-f(x)) = \epsilon>0$ and $B_\rho(S,g)\neq X$. Since $\overline{B_{\rho}(S,f)} \subseteq B_\rho(S,g)$, by $(ii)$, $B_\rho(S,f) \in \mathcal{B}_d(X)$. Let $S' = B_\rho(S,f)$. Since $Id:(X,d) \to (X, \rho)$ is strongly uniformly continuous on $\mathcal{S}_{\mathcal{Z}^+}^\rho$, $\exists$ $\delta > 0$ $B_d(S', \delta) \subseteq B_\rho(S',\epsilon) \subseteq B_\rho(S,g)$. So we have $$B_\rho(S,f) \subseteq B_d(S', \frac{\delta}{3}) \subseteq B_d(S', \frac{\delta}{2})\subseteq B_\rho(S,g).$$  Hence by Theorem \ref{Tsddiffbornologiesmetrics}, the implication $(i)$ follows.
\end{proof} 	

\begin{corollary}\label{AwWijsmandifferentmetrics}
	Let $X$ be a metrizable space. Suppose $d, \rho$ are two compatible metrics on $X$ such  that $(X, \rho)$ is an a.c.m.s. Consider the following statements.
	\begin{enumerate}[(i)]
		\item $\tau_{{AW}_d}^+ \supseteq \tau_{W_\rho}^+$ on $CL(X)$.
		\item each proper $\rho$-closed ball is $d$-bounded.
	\end{enumerate}
	Then $(i) \Rightarrow (ii)$ holds. In addition in $(ii)$, if $Id:(X,d)\to (X, \rho)$ is strongly uniformly continuous on each proper $\rho$-ball, then the converse also holds.  
\end{corollary}
\begin{remark}
	Suppose $(X, \| \cdot \|)$ is a normed linear space, then we have $\sup\{\tau_{W_\rho}: \rho \text{ is an equivalent norm on }X\}$ is equal to \textit{Joly-slice convergence} (see, Theorem 2.4.5 of \cite{ToCCoS} and \cite{joly1973famille}). Since every normed linear space is an a.c.m.s., by Corollary \ref{AwWijsmandifferentmetrics} we get $\tau_{AW_{\| \cdot \|}}$ is finer than the \textit{Joly-slice convergence} as proved in Proposition $3.1.5$ of \cite{ToCCoS}. 
\end{remark}
The next result is a generalized version of Theorem \ref{Tsdunirformlyequimetric}. For two bornologies $\mathcal{S}$, $\mathcal{B}$ on $(X,d)$, we examine when $\tau_{\mathcal{S},d}^+$ is infima of the collection $\{\tau_{\mathcal{B},\rho}^+: \rho \text{ is uniformly equivalent to }d\}$. 	A slight modification of the Theorem \ref{Tsdunirformlyequimetric} yields the proof of this result.
\begin{theorem}\label{Tsbdrho}
	Let $(X,d)$ be a metric space and let $\mathcal{S}, \mathcal{B}$ be two bornologies on $X$. Consider the following statements.
	\begin{enumerate}[(i)]
		\item $\tau_{\mathcal{S},d}^+ \subseteq \tau_{\mathcal{B},\rho}^+$ on $CL(X)$ for each uniformly equivalent metric $\rho$ on $(X,d)$.
		\item For $S \in \mathcal{S}$ and $f,g \in \mathcal{Z}^+$ whenever $\inf_{x \in S}(g(x)-f(x))>0$ and $B_d(S,g) \neq X$, then $B_d(S,f)$ is weakly $\mathcal{B}$-totally bounded with respect to $d$. 
		\end{enumerate}
		Then $(i) \Rightarrow (ii)$ holds. If in $(ii)$, $B_d(S,f)$ is $\mathcal{B}$-totally bounded with respect to $d$, then we also have $(ii) \Rightarrow (i)$. 
\end{theorem}

In \cite{costantini1993metrics}, it is proved that a metrizable space is locally separable if and only if there exists a metric $d$ on $X$ such that $\tau_{W_d}$  is the minimum of the family $\{\tau_{W_\rho}: \rho \text{ is an uniformly equivalent metric to }d\}$. We now prove a similar result for the Attouch-Wets convergence. 
\begin{corollary}\label{AwWijsmancorollary1}
	Let $(X, d)$ be a metric space. Then the following statements are equivalent:
	\begin{enumerate}[(i)]
		\item $\tau_{{AW}_d}^+ \subseteq \tau_{W_\rho}^+$ on $CL(X)$ for each uniformly equivalent metric $\rho$ on $X$;
		\item each proper closed and $d$-bounded set is $d$-totally bounded.
	\end{enumerate}
\end{corollary}
\begin{proof}
	$(i)\Rightarrow (ii)$. Let $B$ be a proper closed and bounded subset of $(X,d)$. Then there exists $y \in X$ and $r > 0$ such that $d(y, B) = 3r>0$. So $B_d(B,2r) \neq X$. By Theorem \ref{Tsbdrho}, $B_d(B,r) \in \mathcal{TB}_d(X)$, thus $B \in \mathcal{TB}_d(X)$. 
	
	$(ii) \Rightarrow (i)$. Let $B \in \mathcal{B}_d(X)$ and $f,g \in \mathcal{Z}^+$ such that $\inf_{x \in B}(g(x)-f(x))>0$ and $B_d(B,g) \neq X$. By Lemma \ref{Bd(X)}, $B_d(B,g) \in \mathcal{B}_d(X)$. Then by $(ii)$, $\overline{B_{d}(B,f)}$ is $d$-totally bounded. By Theorem \ref{Tsbdrho}, $\tau_{{AW}_d}^+ \subseteq \tau_{W_\rho}^+$ on $CL(X)$ for each uniformly equivalent metric $\rho$ on $X$.   
\end{proof}

 In the next corollary we characterize separability of a metrizable space in terms of infima of Wijsman topologies over all uniformly equivalent metrics.  
\begin{corollary}\label{AwWijsmancorollary2}
	Let $X$ be a metrizable space. Then the following statement are equivalent:
	\begin{enumerate}[(i)]
		\item there exists a compatible metric $d$ such that $\tau_{{AW}_d}^+$ is the infimum of the family   $\{\tau_{W_\rho}^+: \rho \text{ is uniformly equivalent metric to }d\}$; 
		 \item $X$ is separable. 
	\end{enumerate}
\end{corollary}

\begin{proof}
	$(i) \Rightarrow (ii)$. Suppose $(i)$ holds. Then by Corollary \ref{AwWijsmancorollary1}, each proper closed and $d$-bounded set is $d$-totally bounded. We show that $X$ is separable. Let $x_0 \in X$. For each $n, m \in \mathbb{N}$, choose $F_{n,m} \in \mathcal{F}(X)$ such that $B_d(x_0, n)\setminus \{x_0\} \subseteq B_d(F_{n,m}, \frac{1}{m})$. Set $D_n = \cup_{m \in \mathbb{N}}F_{n,m}$ and $D = \cup_{n \in \mathbb{N}}D_n \cup \{x_0\}$. It is easy to verify that $D$ is dense in $(X,d)$. Hence $X$ is separable.
	
	$(iii) \Rightarrow (i)$. Suppose $X$ is a separable metrizable space. It is well-known fact that there exists a compatible totally bounded metric $d$ on $X$. So by Corollary \ref{AwWijsmancorollary1}, $\tau_{{AW}_d}^+ = \inf\{\tau_{W_\rho}^+: \rho \text{ is uniformly equivalent to }d\}$ on $CL(X)$.     
\end{proof}

  	  \bibliographystyle{plain}
 	  
 	  \nocite{willard}
 	  \bibliography{reference_file}

\end{document}